\documentclass{amsart}
\usepackage{latexsym}
\usepackage[all]{xypic}
\usepackage{isolatin1}
\numberwithin{equation}{section}

\newtheorem{Lem}{Lemma}[section]
\newtheorem{Prop}[Lem]{Proposition}
\newtheorem{Cor}[Lem]{Corollary}
\newtheorem{Thm}[Lem]{Theorem}

\theoremstyle{definition}

\theoremstyle{remark}
\newtheorem{Rem}[Lem]{Remark}

\renewcommand\o{\otimes}
\DeclareMathOperator\Hom{\operatorname{Hom}}

\DeclareMathOperator\ev{\operatorname{ev}}
\DeclareMathOperator\db{\operatorname{db}}

\newcommand\nt{\diamondsuit}
\newcommand\HMod[4]{{^{#1}_{#3}\mathcal M^{#2}_{#4}}}
\newcommand\LMod[1]{{_{#1}\mathcal M}}
\newcommand\LComod[1]{{^{#1}\mathcal M}}

\newcommand\BiMod[1]{{_{#1}\mathcal M_{#1}}}

\newcommand{\leer}{\operatorname{--}}
\newcommand{\ou}[1]{\mathrel{\mathop{\otimes}_{#1}}}
\newcommand{\co}[1]{\mathrel{\mathop{\Box}_{#1}}}

\newcommand\sw[1]{{}_{(#1)}}
\newcommand\swm[1]{{}_{(-#1)}}
\newcommand\so[1]{^{(#1)}}
\newcommand\som[1]{^{(-#1)}}

\newcommand\ol{\overline}

\newcommand\inv{^{-1}}

\renewcommand\epsilon\varepsilon

\makeatletter
\def\namelabel#1#2{\@bsphack
  \protected@write\@auxout{}%
         {\string\newlabel{#1.nme}{{#2}{#2}}}%
  \@esphack}
\def\nmlabel#1#2{\label{#2}\namelabel{#2}{#1}}
\newcommand\nmref[1]{\ref{#1.nme}\ \ref{#1}}
\makeatother

\begin{document}
\title{Quotients of finite quasi-Hopf algebras}
\author{Peter Schauenburg}
\address{Mathematisches Institut der Universit\"at M\"unchen, 
Theresienstr.~39, 80333~M\"unchen, Germany}
\email{schauen@rz.mathematik.uni-muenchen.de}
\subjclass{16W30}
\keywords{Quasi-Hopf algebra, Nichols-Zoeller theorem, Hopf module}
\begin{abstract}
Let $H$ be a finite-dimensional quasi-Hopf algebra. We show for each
quotient quasibialgebra $Q$ of $H$ that $Q$ is a quasi-Hopf algebra
whose dimension divides the dimension of $H$.
\end{abstract}
\maketitle
\section{Introduction}
In \cite{NicZoe:HAFT} Nichols and Zoeller prove what is now known
as the Nichols-Zoeller Theorem: A finite dimensional Hopf algebra
$H$ over a field $k$ is a free module over every Hopf subalgebra
$K\subset H$. This answers affirmatively one of Kaplansky's conjectures
on Hopf algebras in the finite-dimensional case. The Nichols-Zoeller
Theorem and some related results are an important tool in the 
study of finite-dimensional Hopf algebras.

Quasi-Hopf algebras, introduced by Drinfeld \cite{Dri:QHA}, are a 
generalization of ordinary Hopf algebras that can be motivated most
easily by looking at representation categories: The category of modules
over a Hopf algebra is a monoidal category, with the module
structure on the tensor product over the base field of two modules
given by the diagonal action via comultiplication. The same thing is
still true for quasibialgebras; the difference is that now the
tensor product of representations is associative with an associativity
isomorphism that differs from the ordinary one for vector spaces.

In \cite{Sch:QHAFT} we have proved the direct generalization of 
the Nichols-Zoeller Theorem to quasi-Hopf algebras: A finite-dimensional
quasi-Hopf algebra $H$ over a field $k$ is a free module over every
quasi-Hopf subalgebra $K\subset H$.  The generalization is made 
possible by the introduction of Hopf modules over quasibialgebras
by Hausser and Nill \cite{HauNil:ITQHA}. Hopf modules and the structure
theorem for Hopf modules are a key ingredient in the proof of the
Nichols-Zoeller Theorem (and the subject of a generalization of the
statement of the theorem). Although quasi-Hopf algebras are not
coassociative coalgebras, and thus comodules over them are not 
immediately defined, one can still define Hopf bimodules, and
prove a structure theorem for them.

One of the standard applications of the Nichols-Zoeller theorem is
to investigate the structure of (semisimple) Hopf algebras by
dimension counting: The dimension of a Hopf subalgebra has to 
divide the dimension of the large Hopf algebra, a version of the
classical Lagrange theorem for finite groups. This can serve to narrow
down the possible examples in classification attempts. To be yet
more specific, one standard argument is that the number of one-dimensional
representations of a Hopf algebra $H$ has to divide the dimenision 
of $H$. Unfortunately, the version of the Nichols-Zoeller Theorem
for quasi-Hopf algebras provided in \cite{Sch:QHAFT} does not 
help at all in this situation: The one-dimensional representations
of $H$ are not related to a quasi-Hopf subalgebra, but rather
to a quotient quasi-Hopf algebra. For a quotient Hopf algebra $Q$
of an ordinary Hopf algebra $H$, the Nichols-Zoeller Theorem 
implies by duality that $H$ is a cofree $Q$-comodule. 
In the quasi-Hopf
case, this does not even make sense to ask, since $Q$ is not
a coassociative coalgebra, hence $H$ is not a comodule in the usual
sense.

We shall nevertheless prove that the dimension of a finite-dimensional
quasi-Hopf algebra $H$ is divisible by the dimension of any
quotient quasi-Hopf algebra $Q$. The key to this is the construction
of an inclusion of quasi-Hopf algebras with the same ratio of 
dimensions as that between $H$ and $Q$. In the ordinary Hopf case
this is found by simply dualizing. In the quasi-Hopf case, we will
find that $H$ and $Q^*$ generate a quasi-Hopf subalgebra 
$D(Q;H)$ in the
Drinfeld double $D(H)$ of $H$, and $\dim D(Q;H)=\dim H\dim Q$,
while $\dim D(H)=(\dim H)^2$.

Without doubt, one could verify the claims just made on $D(Q;H)$
by direct calculations with the rather complicated quasi-Hopf algebra
structure of the Drinfeld double given by Hausser and Nill \cite{HauNil:DQQG}
(while the earlier description by Majid \cite{Maj:QDQHA} is perhaps
to indirect for this purpose). Instead, we will do a closer analysis
of $D(Q;H)$ in \nmref{sec:reldo}, giving parallel interpretations for
its modules to those of the double $D(H)$. This will allow us to 
show our claims without calculating much. Moreover, we will be able
in \nmref{sec:lagrange} to show more than the ``Lagrange'' statement
that $\dim Q$ divides $\dim H$: We will show that $\dim Q$ divides
the dimension of any Hopf module in $\HMod Q{}HH$, parallel to the
results of Nichols and Zoeller who also show a freeness result for
Hopf modules, not only for the Hopf algebras themselves.

In addition to the results summarized so far, we will show in 
\nmref{sec:ant} that any quotient quasibialgebra 
of a finite-dimensional 
quasi-Hopf algebra is a quotient quasi-Hopf algebra itself. In
the case of ordinary Hopf algebras, this is due to Nichols
\cite{Nic:QHA}, whose arguments we will vary in the necessary
manner, replacing a canonical map $H\o H\rightarrow H\o H$
by its quasi-Hopf version due to Drinfeld. Nichols' result
applies (by duality) to subbialgebras of finite-dimensional 
Hopf algebras, whereas it is not true for subquasibialgebras
of finite-dimensional quasi-Hopf algebras, although we will
give a positive result under additional hypotheses.

In an appendix, we will give a categorical proof, rather free of 
computations, of the canonical isomorphism $H\o H\rightarrow H\o H$
given by Drinfeld in \cite{Dri:QHA} and used crucially in 
\nmref{sec:ant}.
\section{Antipodes for quotients and subobjects}\nmlabel{Section}{sec:ant}
Recall that 
a quasibialgebra $H=(H,\Delta,\epsilon,\phi)$ consists of an algebra $H$,
algebra maps $\Delta\colon H\rightarrow H\o H$ and $\epsilon\colon H\rightarrow k$,
and an invertible element $\phi\in H^{\o 3}$, the associator, such that
\begin{gather}
(\epsilon\o H)\Delta(h)=h=(H\o\epsilon)\Delta(h),\\
(H\o\Delta)\Delta(h)\cdot\phi=\phi\cdot(\Delta\o H)\Delta(h)\label{quasicoass},\\
(H\o H\o\Delta)(\phi)\cdot(\Delta\o H\o H)(\phi)=(1\o\phi)\cdot(H\o\Delta\o H)(\phi)\cdot(\phi\o 1)\label{cocycle},\\
(H\o\epsilon\o H)(\phi)=1
\end{gather}
hold for all $h\in H$.
We will write $\Delta(h)=:h\sw 1\o h\sw 2$, 
$\phi=\phi\so 1\o\phi\so 2\o\phi\so 3$, and
$\phi\inv=\phi\som 1\o\phi\som 2\o\phi\som 3$. 

We define a morphism or map of quasibialgebras from 
$(H,\phi)$ to $(L,\psi)$ to be an algebra map $f\colon H\rightarrow L$
compatible with comultiplication $\Delta f=(f\o f)\Delta$ and with
the counit $\epsilon f=\epsilon$, and finally satisfying
$(f\o f\o f)(\phi)=\psi$.

A subquasibialgebra is a subalgebra $L$ of $H$ that has a 
quasibialgebra structure for which the inclusion is a quasibialgebra
map. Equivalently $L$ is a subalgebra of $H$ satisfying
$\Delta(H)\subset H\o H$ and $\phi\in H\o H\o H$. A quotient
quasibialgebra of $H$ is a quasibialgebra $Q$ for which there is
a surjective quasibialgebra map $\nu\colon H\rightarrow Q$. Equivalently,
$Q$ is isomorphic to $H/I$ for an ideal $I\subset H$ satisfying
$\Delta(I)\subset I\o H+H\o I$ and $\epsilon(I)=0$; a coassociator
for $Q=H/I$ is then the canonical image of $\phi$ in $Q\o Q\o Q$.

Recall further that a quasi-antipode for a quasibialgebra $(H,\phi)$
is a triple $(S,\alpha,\beta)$ in which $\alpha,\beta\in H$,
and $S$ is an anti-algebra endomorphism of $H$
satisfying
\begin{align*}
  S(h\sw 1)\alpha h\sw 2&=\epsilon(h)\alpha, & h\sw 1\beta S(h\sw 2)&=\epsilon(h)\beta,\\
  \phi\so 1\beta S(\phi\so 2)\alpha\phi\so 3&=1,&S(\phi\som 1)\alpha\phi\som 2\beta\phi\som 3&=1
\end{align*}
for $h\in H$. A quasi-Hopf algebra is a quasibialgebra
with a quasi-antipode. Note that our definition differs from 
Drinfeld's in that we do not require the antipode to be bijective.
For finite-dimensional quasibialgebras it was recently shown
by Bulacu and Caenepeel \cite{BulCae:IDQHAA} that bijectivity
of the antipode is automatic.

We define a quasibialgebra map $f\colon H\rightarrow H'$
between quasi-Hopf algebras
$(H,\phi,S,\alpha,\beta)$ and $(H',\phi',S',\alpha',\beta')$ to
be a quasi-Hopf algebra map if $S'f=fS$, $f(\alpha)=\alpha'$,
and $f(\beta)=\beta'$. We define a quasi-Hopf subalgebra to be
a subquasibialgebra  $L\subset H$ that has a quasi-antipode so that
the inclusion is a quasi-Hopf algebra map. Equivalently, $L$
is a subquasibialgebra such that $S(L)\subset L$, and $\alpha,\beta\in L$.
Further we define a quotient quasi-Hopf algebra of $H$ to be a
quasi-Hopf algebra $Q$ with a surjective quasi-Hopf algebra map
$\nu\colon H\rightarrow Q$. Equivalently, a quotient quasi-Hopf
algebra is a quotient quasibialgebra $Q\cong H/I$ such that
$S(I)=I$. Then $S$ induces an antiautomorphism on the quotient, 
which is a quasiantipode together with the canonical images of
$\alpha$ and $\beta $ in $Q$.

\begin{Thm}\nmlabel{Theorem}{quotant}
Let $H$ be a quasi-Hopf algebra over a field $k$, and let
$Q$ be a finite dimensional 
quotient quasibialgebra of $H$. Then $Q$ is a quotient quasi-Hopf
algebra.
\end{Thm}
\begin{proof}
By \cite[Prop.~1.5]{Dri:QHA} we have an isomorphism
\begin{equation}\label{can}
\varphi\colon H\o H\ni g\o h\mapsto gS(\phi\som 1)\alpha\phi\som 2h\sw 1\o \phi\som 3h\sw 2\in H\o H
\end{equation}
with inverse given by
$$\varphi\inv(g\o h)=g\phi\so 1\beta S(h\sw 1\phi\so 2)\o h\sw 2\phi\so 3$$
Let $\nu\colon H\ni h\mapsto \ol h\in Q$ 
denote the
canonical epi.
Define 
$$\ol\varphi\colon Q\o Q\ni p\o q\mapsto p\ol{S(\phi\som 1)\alpha\phi\som 2}q\sw 1\o q\sw 2\in Q\o Q.$$
Then the diagram
$$\xymatrix{H\o H\ar[r]^-{\varphi}\ar[d]^{\nu\o\nu}&H\o H\ar[d]^{\nu\o\nu}\\
  Q\o Q\ar[r]^-{\ol\varphi}&Q\o Q}$$
commutes, and thus $\ol\varphi$ is onto. Since $Q$ is finite dimensional,
$\ol\varphi$ is an isomorphism. From the diagram we conclude that
$$\ol\varphi\inv(\ol g\o\ol h)=(\nu\o\nu)\varphi\inv(g\o h),$$
hence
$$\theta\colon Q\ni q\mapsto(Q\o\epsilon)\ol\beta\inv(1\o q)\in Q$$
satisfies
$\theta(\ol h)=\ol{\beta S(h)}$ for all $h\in H$.
Now define
$$\ol S\colon Q\ni q\mapsto \ol{S(\phi\som 1)\alpha\phi\som 2}\theta(q\ol{\phi\som 3})\in Q.$$
Then for $h\in H$ we have
$$\ol S(\ol h)=\ol{S(\phi\som 1)\alpha\phi\som 2}\cdot\ol{\beta S(h\phi\som 3)}
=\ol{S(\phi\som 1)\alpha\phi\som 2\beta S(\phi\som 3)}\cdot\ol{S(h)}=\ol{S(h)},$$
showing that $S$ maps the kernel of $\nu$ into itself.
\end{proof}
\begin{Rem}
  The direct analog of \nmref{quotant} for subquasibialgebras 
  instead of quotients is false for quite trivial reasons. To see
  this consider a quasi-Hopf algebra $(H,\phi,S,\alpha,\beta)$ and
  a quasi-Hopf subalgebra $K\subset H$. By the remark following
  the definition of a quasi-Hopf algebra in \cite{Dri:QHA}, 
  we can obtain another quasi-Hopf structure
  $(H,\phi,S',\alpha',\beta')$ for any unit $u\in H$ by setting
  $S'(h)=uS(h)u\inv$, $\alpha'=u\alpha$, $\beta'=\beta u\inv$,
  while leaving $\phi$ unchanged.
  Of course, it may happen that $K$ is not a quasi-Hopf subalgebra
  for this new quasi-Hopf structure (for example, if $\alpha$ is 
  a unit, and $u\not\in K$). 
\end{Rem}
  
  However, we can provide quasiantipodes for subquasibialgebras
  under some additional assumptions:
\begin{Prop}\nmlabel{Proposition}{subprop}
  Let $H$ be a quasi-Hopf algebra with 
  coassociator $\phi$ and quasi-antipode $(S,\alpha,\beta)$.

  Let $K\subset H$ be a finite-dimensional subquasibialgebra such that
  $$S(\phi\som 1)\alpha\phi\som 2\o\phi\som 3\in K\o K.$$
  Then $K$ is a quasi-Hopf subalgebra.
\end{Prop}
\begin{proof}
  Again we consider the canonical map $\varphi$ from 
  \eqref{can}. By our extra assumptions, we see that 
  $\varphi(K\o K)\subset K\o K$, and consider the
  map $\varphi'\colon K\o K\rightarrow K\o K$ given by restricting
  $\varphi$. It is injective since $\varphi$ is, hence bijective
  by finite dimensionality. The inverse of $\varphi'$ is given
  by the restriction of $\varphi\inv$, so we see that for 
  $x\in K$
  $$K\ni(K\o\epsilon)\varphi\inv(1\o x)=\beta S(x),$$
  hence in particular $\beta\in K$, and for all $x\in K$
  $$K\ni S(\phi\som 1)\alpha\phi\som 2\beta S(x\phi\som 3)
         =S(x).$$
  Finally $\alpha=(K\o\epsilon)\varphi(1\o 1)\in K$,
  so $K$ is a quasi-Hopf subalgebra.
\end{proof}
\begin{Rem}
  As a special case of \nmref{subprop}, a finite-dimensional subquasibialgebra
  $K\subset H$ of a quasi-Hopf algebra $H$
  is a quasi-Hopf subalgebra provided that it contains
  a subquasibialgebra $L\subset K$ which is a quasi-Hopf subalgebra
  of $H$.
\end{Rem}

\section{The partial double}\nmlabel{Section}{sec:reldo}
Throughout the section, we let $H$ denote a quasi-Hopf algebra.
The key property of a quasibialgebra is that its modules form
a monoidal category: The tensor product of $V,W\in\LMod H$ is
their tensor product $V\o W$ over $k$, endowed with the 
diagonal module structure $h(v\o w)=h\sw 1v\o h\sw 2w$; the 
neutral object is $k$ with the trivial module structure given
by $\epsilon$. The associativity isomorphism in the category
is 
$$(U\o V)\o W\ni u\o v\o w\mapsto\phi\so 1u\o\phi\so 2v\o\phi\so 3w\in U\o(V\o W)$$
for $U,V,W\in\LMod H$.

The opposite of a quasibialgebra and the
tensor product of two quasibialgebras are naturally quasibialgebras.
Thus $\BiMod H$ is also a monoidal category, with 
associativity isomorphism
$$(U\o V)\o W\ni u\o v\o w\mapsto\phi\so 1u\phi\som 1\o\phi\so 2v\phi\som 2\o\phi\so 3w\phi\som 3\in U\o(V\o W).$$

We will make free use of the formalism of (co)algebra and (co)module
theory within monoidal categories. When $C,D$ are coalgebras
in $\BiMod H$, we will use the abbreviations $\HMod C{}HH,\HMod{}DHH,
\HMod CDHH$ for the categories of left $C$-comodules, right $D$-comodules,
and $C$-$D$-bicomodules within the monoidal category $\BiMod H$.

We see that $H$ itself is a coassociative coalgebra within the
monoidal category $\BiMod H$. Thus we can define a
Hopf module $M\in\HMod{}HHH$ to be a right $H$-comodule within
the category $\BiMod H$. Written out explicitly, this definition
is the same as that of Hausser and Nill \cite[Def.3.1]{HauNil:ITQHA}.
Hausser and Nill have also proved a structure theorem for
such Hopf modules, which says that the functor
$$  \mathcal R\colon \LMod H\ni  V\mapsto {_\cdot V}\o{_\cdot H_\cdot^\cdot}
\in\HMod{}HHH$$
is an equivalence of categories. We have used this equivalence
as a basis for a description of the Drinfeld double of $H$
in \cite{Sch:HMDQHA}.
In \cite[Expl.4.10]{Sch:AMCGHSP} we have repeated this description
with a general $\mathcal C$-categorical technique, which we shall
follow once more now to obtain a relative double $D(L;H)$ for any
quasibialgebra map $\nu\colon H\rightarrow L$.

For any right $H$ comodule $M$ in $\BiMod H$ and any $P\in \BiMod H$
we can form the right $H$-comodule $P\o M$ in $\BiMod H$, which
gives us a functor $\BiMod H\times \HMod{}HHH\rightarrow\HMod{}HHH$
that makes $\HMod{}HHH$ into a left $\BiMod H$-category in the sense
of Pareigis \cite{Par:NARMTII}. Being equivalent to $\HMod{}HHH$, 
the category $\LMod H$ is then also a left $\BiMod H$-category,
which means that we have a functor 
$\nt\colon\BiMod H\times\LMod H\rightarrow \LMod H$ and
a coherent natural isomorphism $\Omega\colon(P\o Q)\nt V\rightarrow P\nt(Q\nt V)$
for $P,Q\in\BiMod H$ and $V\in\LMod H$. 

Now let $C$ be a coalgebra in $\BiMod H$. 
Since $\LMod H$ is a left $\BiMod H$-category, it makes sense 
(see \cite{Par:NARMTII}) to
talk about $C$-comodules within $\LMod H$, which form a 
category $^C\left(\LMod H\right)$, which in our situation is 
naturally equivalent to $\HMod CHHH$, with the equivalence
$^C\left(\LMod H\right)\cong \HMod CHHH$ induced by 
$\mathcal R$.

By \cite[Thm.3.3]{Sch:AMCGHSP} and the remarks preceding it,
$\nt$ induces a functor $\BiMod H\ni P\mapsto P\nt H\in\BiMod H$,
and we have an isomorphism
$(P\nt H)\ou HV\cong P\nt V$, natural in $P\in\BiMod H$ and
$V\in \LMod H$.

By \cite[Cor.3.8]{Sch:AMCGHSP}, $C\nt H$ has an $H$-coring
structure in such a way that one has 
a category equivalence
$\LComod{C\nt H}\cong {^C\left(\LMod H\right)}$ that
commutes with the underlying functors to $\LMod H$.
(Here $\LComod{C\nt H}$ denotes the category of left comodules
over the $H$-coring $C\nt H$.) 
For any coalgebra morphism
$f\colon C\rightarrow D$ in $\BiMod H$ we obtain a commutative
diagram of functors
$$\xymatrix{\LComod{C\nt H}\ar[r]\ar[d]^{\LComod{f\nt H}}
    &{^C\left(\LMod H\right)}\ar[d]^{^f\left(\LMod H\right)}\ar[r]^{\mathcal R}
    &\HMod CHHH\ar[d]^{\HMod fHHH}\\
    \LComod{D\nt H}\ar[r]
    &{^D\left(\LMod H\right)}\ar[r]^{\mathcal R}
    &\HMod DHHH}$$
for the $H$-coring map $f\nt H\colon C\nt H\rightarrow D\nt H$.

From \cite[Sec.4]{Sch:HMDQHA} we know that 
$P\nt V\cong P\o V$ as vector spaces, functorially in 
$P\in\BiMod H$ and $V\in\LMod H$. In particular
we have $C\nt H\cong C\o H$ as right $H$-modules, functorially in 
the coalgebra $C$ in $\BiMod H$.

Now assume that $C$ is finite dimensional. Then $C\nt H$ is a 
finitely generated free right $H$-module, so the dual
algebra $\Hom_{-H}(C\nt H,H)$ of the $H$-coring $C\nt H$ 
fulfills $\LMod{(C\nt H)^\vee}\cong\LComod{C\nt H}$. Note
that $(C\nt H)^\vee\cong H\o C^*$ as vector spaces, functorially
in $C$.

Finally, we specialize to a class of examples of coalgebras
in $\BiMod H$. Whenever $\nu\colon H\rightarrow L$ is a morphism
of quasibialgebras, we can consider $L$ as a coalgebra in $\BiMod H$
with respect to its usual comultiplication, and the $H$-bimodule
structure induced along $\nu$. We write $D(L;H):=(L\nt H)^\vee$ for
the dual algebra of the coring $L\nt H$. 
Note that $D(L;H)\cong H\o L^*$ as vector spaces. The isomorphism
is natural in $L$, meaning that for any morphism $f\colon L\rightarrow M$
of quasibialgebras, the induced morphism
$D(f;H)\colon D(M;H)\rightarrow D(L;H)$ corresponds to 
$H\o f^*\colon H\o M^*\rightarrow H\o L^*$. In particular
it is surjective (resp.\ injective) if $f$ is injective 
(resp.\ surjective). 

The same calculations as those made to prove 
\cite[Lem.3.2]{Sch:HMDQHA} prove more generally that 
$\HMod LHHH$ is a monoidal category, the tensor product of
$M,N\in\HMod LHHH$ being $M\ou HN$ with the left and right 
comodule structures
\begin{gather*}
  M\ou HN\ni m\o n\mapsto m\swm 1n\swm 1\o m\sw 0\o n\sw 0\in L\o (M\ou HN)\\
  M\ou HN\ni m\o n\mapsto m\sw 0\o n\sw 0\o m\sw 1n\sw 1\in (M\ou HN)\o H.
\end{gather*}
Note that the underlying functor $\HMod LHHH\rightarrow \HMod {}HHH$
is monoidal.

The equivalence $\mathcal R$ is a monoidal equivalence. 
We make $\LMod{D(L;H)}$ a monoidal category in such a way that
the equivalence $\LMod{D(L;H)}\cong\HMod LHHH$ is a monoidal functor.
Thus, for any quasibialgebra maps $H\xrightarrow\nu L\xrightarrow fM$
we obtain a commutative diagram of monoidal functors
$$\xymatrix{\LMod{(M\nt H)^\vee}\ar[rr]^-{\LMod{(f\nt H)^\vee}}\ar[d]
    &&\LMod{(L\nt H)^\vee}\ar[rr]^-{\LMod{(\epsilon\nt H)^\vee}}\ar[d]
    &&\LMod H\ar[d]\\
    \HMod MHHH\ar[rr]^{\HMod fHHH}&&\HMod LHHH\ar[rr]&&\HMod {}HHH}
$$
By a trivial modification of the monoidal category structures
(cf.\ \cite[Rem.5.3]{Sch:BNRSTHB} for an analogous trick) 
we can make sure that the functors in 
the top row are strict monoidal functors. This implies that
$D(L;H)$ is a quasibialgebra, and that 
$D(f;H)\colon D(M;H)\rightarrow D(L;H)$ is a quasibialgebra map.

\section{Lagrange's theorem}\nmlabel{Section}{sec:lagrange}
Consider a finite-dimensional quasi-Hopf algebra $H$, and a quotient
quasi-Hopf algebra $Q$. 
As a particular case of the constructions in the preceding section,
we obtain injective quasibialgebra maps $H\rightarrow D(Q;H)\rightarrow D(H;H)$.
Note that $D(H;H)=D(H)$ is the Drinfeld double of $H$, which is a
quasi-Hopf algebra. It follows from \nmref{subprop}
that $D(Q;H)$ is a quasi-Hopf algebra
as well, and a quasi-Hopf subalgebra of $D(H)$. 
By \cite{Sch:QHAFT}, $D(H)$ is a free $D(Q;H)$-module, so in 
  particular $\dim D(Q;H)=\dim Q\dim H$ divides
  $\dim D(H)=(\dim H)^2$. Cancelling $\dim H$ we get:
\begin{Cor}
  Let $H$ be a finite-dimensional quasi-Hopf algebra, and $Q$
  a quotient quasibialgebra of $H$. Then $\dim Q$ divides $\dim H$.
\end{Cor}
As an immediate application we have:
\begin{Cor}
  Let $H$ be a finite-dimensional quasi-Hopf algebra. Then the
  number of one-dimensional representations of $H$ divides
  $\dim H$.
\end{Cor}
\begin{proof}
  We pass to the dual picture: Representations of $H$ are comodules
  over the dual coalgebra $H^*$. One-dimensional comodules correspond
  to grouplike elements of $H^*$. These grouplikes span a 
  sub-coquasibialgebra of $H^*$, which corresponds to a quotient
  quasibialgebra of $H$, whose dimension divides the dimension of $H$
  by the preceding corollary.
\end{proof}
Nichols and Zoeller \cite{NicZoe:HAFT} do not only prove a 
freeness theorem for Hopf algebra inclusions $K\subset H$, but
also a freeness theorem for Hopf modules in $\HMod H{}K{}$.
In particular, for any finite-dimensional 
$M\in\HMod H{}K{}$,
they show that $\dim K|\dim M$.
Suppose that $H$ is a finite-dimensional quasi-Hopf algebra and
$Q$ a quotient quasibialgebra. Then $Q$ is a coalgebra in the
monoidal category $\BiMod H$, so that we have a well-defined notion
of Hopf module in $\HMod Q{}HH$ (while $\HMod Q{}H{}$ is not defined).
\begin{Prop}
  Let $H$ be a finite-dimensional quasi-Hopf algebra, $Q$ a 
  quotient quasibialgebra of $H$, and 
  $M\in\HMod Q{}HH$. Then $\dim Q$ divides $\dim M$.
\end{Prop}
\begin{proof}
Consider the commutative diagram of functors (first without the dotted arrows)
$$\xymatrix{\LMod{D(H)}\ar[r]^-{\mathcal R}\ar@<-.5ex>[d]_{\mathcal V}&\HMod HHHH\ar@<-.5ex>[d]_{\mathcal U}\\
      \LMod{D(Q;H)}\ar@<-.5ex>@{..>}[u]\ar[r]^-{\mathcal R}&\HMod QHHH\ar@<-.5ex>@{..>}[u]}$$
in which the functor $\mathcal V$ is the underlying functor induced
by the inclusion $D(Q;H)\rightarrow D(H)$, and $\mathcal U$
is the underlying functor induced by the projection $H\rightarrow Q$.
The horizontal arrows are induced by the category equivalence
$\mathcal R\colon\LMod H\rightarrow \HMod{}HHH$. Since they are
equivalences, the diagram also commutes for the dotted arrows, if these
denote the right adjoint functors to $\mathcal V$ and $\mathcal U$.
Now the right adjoint to the underlying functor $\mathcal U$
is given by cotensor product with $H$, taken within the monoidal
category $\BiMod H$ (dually to the induction functor for an 
algebra inclusion), whereas the right adjoint to $\mathcal V$
is the usual coinduction functor for the algebra inclusion
$D(Q;H)\subset D(H)$. 
This means that for $W\in\LMod{D(Q;H)}$ we have
$$H\co Q\mathcal R(W)\cong\mathcal R(\Hom_{D(Q;H)}(D(H),W)).$$
In particular, since $D(H)$ is a free $D(Q;H)$-module
of rank $\dim H/\dim Q$, and since $\mathcal R$ multiplies 
dimensions by $\dim H$, we have
$$\dim(H\co Q M)=\frac{\dim H}{\dim Q}\dim M$$
whenever $M\in\HMod QHHH$ is finite dimensional.

Now consider a finite-dimensional $V\in\HMod Q{}HH$. Then 
we can tensor $V$ with $H\in\HMod{}HHH$ to obtain
$M=V\o H\in\HMod QHHH$. Calculating within the monoidal category
$\BiMod H$ we have
$$(H\co QV)\o H\cong H\co Q(V\o H)=H\co QM.$$
It follows that 
$$\dim H\dim (H\co QV)=\frac{\dim H}{\dim Q}\dim M$$
or $\dim (H\co QV)=\frac{\dim H}{\dim Q}\dim V$.

But on the other hand $H\co QV\in\HMod H{}HH$. By the left-right
switched version of the structure theorem of 
Hausser and Nill, any Hopf module in $\HMod H{}HH$ is a free
left $H$-module, so that
$\dim H$ divides $\dim(H\co QV)$. Thus
$\dim Q$ divides $\dim V$.
\end{proof}

\appendix
\section{A dogmatic proof of a formula of Drinfeld}
In this section we return to the isomorphism $\varphi$ from 
equation \eqref{can}. Its proof in Drinfeld's paper \cite{Dri:QHA}
is not particularly hard, but does involve a calculation with 
the coassociator element of $H$, 
the pentagon identity \eqref{cocycle} and the identities defining a
quasi-antipode. The ``dogma'' alluded to in this section's title
says that such calculations should be banned. After all, the pentagon
identity precisely ensures that $\LMod H$ is a monoidal category,
and the quasi-antipode axioms precisely ensure that the dual 
vector space of a finite-dimensional module can be made into a 
dual object inside that monoidal category. Since the axioms
are more or less equivalent to the categorical properties, no further
reference to the former should be necessary, and using the 
latter should lead to easier and more conceptual proofs.

So now we let $H$ be a finite-dimensional quasi-Hopf algebra
and $V\in \LMod H$ a finite-dimensional $H$-module. 
Then $V$ has the dual object $V^*\in\LMod H$, with module
structure given by the transpose of the right module structure
$V_S$ induced via the quasi-antipode, evaluation
$\ev\colon V^*\o V\ni\varphi\o v\mapsto \varphi(\alpha v)\in k$
and dual basis $\db\colon k\rightarrow V\o V^*$,
$\db(1)=\beta v_i\o v^i$, where $v_i$ and $v^i$ are a pair of dual
bases in $V$ and $V^*$, and summation is suppressed.
It follows that 
the functor
$$\mathcal F\colon \LMod H\ni W\mapsto {_\cdot W}\o {_\cdot V}\in\LMod H$$
has the right adjoint
$$\mathcal G\colon\LMod H\ni X\mapsto X\o V^*\in\LMod H$$
with unit $u$ and counit $c$ of the adjunction given by
\begin{gather*}
u=\left(X\xrightarrow{X\o\db}X\o(V\o V^*)\xrightarrow{\Phi\inv}(X\o V)\o V^*\right)\\
c=\left((X\o V^*)\o V\xrightarrow{\Phi}X\o(V^*\o V)\xrightarrow{X\o\ev}X\right).
\end{gather*}

We have the standard isomorphism
$$\Gamma=\left( X\o V^*\overset\cong\rightarrow \Hom(V,X)\overset\cong\rightarrow\Hom_{H-}({_\cdot H}\o V,X)\right)$$
with $\Gamma(x\o\varphi)(h\o v)=h\varphi(v)x$ and
$\Gamma\inv(f)=f(1\o v_i)\o v^i$. 
Via the isomorphism $\Gamma$ we have another right adjoint 
$$\mathcal G'\colon\LMod H\ni X\mapsto \Hom_{H-}({_\cdot H}\o V,X)\in\LMod H.$$
The unit $u'$ and counit $c'$ of the adjunction are given by
\begin{gather*}
  u'=\left(X\xrightarrow{u}(X\o V)\o V^*\xrightarrow{\Gamma}\Hom_{H-}(Q,X\o V)\right)\\
  c'=\left(\Hom_{H-}(Q,X)\o V\xrightarrow{\Gamma\inv\o V}(X\o V^*)\o V\xrightarrow{c}X\right)
\end{gather*}
so that
\begin{multline*}
  u'(x)(h\o v)=\Gamma(u(x))(h\o v)
    =\Gamma(\phi\som 1x\o\phi\som 2\beta v_i\o\phi\som 3v^i)(h\o v)
    \\=h(\phi\som 3v^i)(v)\cdot(\phi\som 1 x\o\phi\som 2\beta v_i)
    =hv^i(v)\cdot(\phi\som 1x\o\phi\som 2\beta S(\phi\som 3)v_i)
    \\h\sw 1\phi\som 1x\o h\sw 2\phi\som 2\beta S(\phi\som 3)v
\end{multline*}
and
\begin{multline*}
  c'(f\o v)=c(f(1\o v_i)\o v^i\o v)
    =\phi\so 1f(1\o v_i)(\phi\so 2v^i)(\alpha\phi\so 3v)
    \\=f(\phi\so 1\o v_i)v^i(S(\phi\so 2)\alpha\phi\so 3v)
    =f(\phi\so 1\o S(\phi\so 2)\alpha\phi\so 3v).
\end{multline*}
One checks that the relevant $H$-module structure on $\Hom_{H-}(H\o V,X)$
(making  $\Gamma$ an $H$-module map)
is given by 
$(hf)(g\o v)=f(gh\sw 1\o S(h\sw 2)v)$. In other words, we have
$\mathcal G'(X)=\Hom_{H-}(Q,X)$ for the 
$H$-bimodule $Q={_\cdot H_\cdot}\o {_\cdot(V_S)}$.
In particular $\mathcal G'$ is the right adjoint in the standard
hom-tensor adjunction with left adjoint $\mathcal F'=Q\ou H(\leer)$.
We denote the counit and unit of that standard adjunction by
$c'',u''$. Now left adjoints are unique, so that we get 
mutually inverse isomorphisms
\begin{gather*}
\Lambda=\left(Q\ou HV\xrightarrow{Q\ou Hu'}Q\ou H\Hom_{H-}(Q,X\o V)\xrightarrow{c''}X\o V\right)\\
\Lambda\inv=\left(X\o V\xrightarrow{u''\o V}\Hom_{H-}(Q,Q\ou HX)\o V\xrightarrow{c'}Q\ou HX\right).
\end{gather*}
We compute
$$\Lambda(h\o v\o x)=u'(x)(h\o v)
=h\sw 1\phi\som 1x \o h\sw 2\phi\som 2\beta S(\phi\som 3)v$$
and
\begin{multline*}
  \Lambda\inv(x\o v)
    =c'(u''(x)\o v)
    =u''(x)(\phi\so 1\o S(\phi\so 2)\alpha\phi\so 3v)
    =\phi\so 1\o S(\phi\so 2)\alpha\phi\so 3v\o x
\end{multline*}
Finally, we specialize $V=X=H$, and identify $Q\ou HH\cong Q=H\o H$
to find
$\Lambda(h\o g)=h\sw 1\phi\som 1\o h\sw 2\phi\som 2\beta S(\phi\som 3)g$
and 
$\Lambda\inv (h\o g)
=\phi\so 1h\sw 1\o S(h\sw 2)S(\phi\so 2)\alpha\phi\so 3g$.

We have obtained the version of \eqref{can} for the opposite and
coopposite quasi-Hopf algebra to $H$.

\end{document}